\documentclass[11pt,letterpaper]{article}
\usepackage{ifthen,latexsym,amssymb,amsmath}

\newcommand{\I}[1]{{\mathbb #1}}

\newcommand{\e}{\varepsilon}
\newcommand{\floor}[1]{\lfloor #1\rfloor}
\newcommand{\me}{{\mathrm e}}
\renewcommand{\mid}{:}

\newif\ifnotesw\noteswtrue

\newcommand{\beq}[1]{\begin{equation}\label{eq:#1}}
\newcommand{\eeq}{\end{equation}}
\newcommand{\req}[1]{\textrm{(\ref{eq:#1})}}

\newtheorem{theorem}{Theorem}
\newcommand{\bth}[2][nothing]{\ifthenelse{\equal{#1}{nothing}}
 {\begin{theorem}} {\begin{theorem}[#1]}\label{th:#2}}

\newtheorem{lemma}[theorem]{Lemma}
\newcommand{\blm}[2][nothing]{\ifthenelse{\equal{#1}{nothing}}
 {\begin{lemma}} {\begin{lemma}[#1]}\label{lm:#2}}

\newtheorem{problem}[theorem]{Problem}
\newcommand{\bpr}[2][nothing]{\ifthenelse{\equal{#1}{nothing}}
 {\begin{problem}} {\begin{problem}[#1]}\label{pr:#2}}

\newcommand{\bpf}[1][Proof.]{\smallskip\noindent{\it #1} }
\newcommand{\qed}{\nolinebreak\mbox{\hspace{5 true pt}%
  \rule[-0.85 true pt]{3.9 true pt}{8.1 true pt}}}

\newcommand{\epf}{\qed \medskip}

\newcommand{\brm}{\smallskip\noindent{\bf Remark.} }

\newcommand{\OPdata}{Oleg Pikhurko\\
Department of Mathematical Sciences\\
Carnegie Mellon University\\
Pittsburgh, PA 15213-3890\\
Web: {\tt http://www.math.cmu.edu/\symbol{126}pikhurko}}

\begin{document}

\newcommand{\ex}{\mbox{\rm ex}}

\author{\OPdata}
\title{Remarks on a Paper by Y.~Caro and R.~Yuster
on Tur\'an Problem\\ (Revised)}
\maketitle

\begin{abstract}
 For a graph $F$ and a function $f:\I N\to\I R$, let
$e_f(F)=\sum_{x\in V(F)} f(d(x))$ and let $\ex_f(n,F)$ be the maximum
of $e_f(G)$ over all $F$-free graphs $G$ with $n$ vertices.

Suppose that $f$ is a non-decreasing function with the property that
for any $\e>0$ there is $\delta>0$ such that for any $n\le m\le
(1+\delta)n$ we have $f(m)\le (1+\e)f(n)$. Under this assumption we
prove that the asymptotics of $\ex_f(n,F)$, where $F$ is a fixed
non-bipartite graph and $n$ tends to the infinity, can be computed by
considering complete $(\chi(F)-1)$-partite graphs only.

This research was motivated by a paper of Y.~Caro and R.~Yuster
[\emph{Electronic J.\ Combin.}, \textbf{7} (2000)] who studied the
case when $f:x\mapsto x^\mu$ is a power function.\end{abstract}

\section{Introduction}

Let $\I N$ denote the set of non-negative integers and $\I R$ the set
of reals. Let $f:\I N\to\I R$ be an arbitrary function.

For a graph $F$ define $e_f(F)=\sum_{x\in V(F)} f(d(x))$, where $d(x)$
denotes the degree of a vertex $x$. For example, for $f:x\mapsto \frac
x2$ we have $e_f(F)=e(F)$; thus $e_f(F)$ can be viewed as a
generalization of the size of $F$.

Define $\ex_f(n,F)$ to be the maximal value of $e_f(G)$ over all
$F$-free graphs $G$ of order $n$. This mimics the definition of the
usual Tur\'an function $\ex(n,F)$. The special case when $f$ is the
\emph{power function} $P_\mu:x\mapsto x^\mu$, with integer $\mu\ge 1$,
appears in the paper of Caro and Yuster~\cite{caro+yuster:00} which
was the main motivation for the present research.

Let $\ex_f'(n,F)$ be the maximum of $e_f(H)$ over all complete
$(\chi(F)-1)$-partite graphs of order $n$. Clearly, we have
 \beq{s1}
 \ex_f'(n,F)\le \ex_f(n,F).
 \eeq

Caro and Yuster~\cite{caro+yuster:00} proved that
$\ex_{P_\mu}(n,K_r)=\ex_{P_\mu}'(n,K_r)$ for $r\ge 3$. However, they
incorrectly claimed that the \emph{Tur\'an} graph $T_{r-1}(n)$, which
has part sizes almost equal, is always optimal. (A revised
version~\cite{caro+yuster:04} of their paper appeared.) This mistake
also appears in the first version of the current paper.

First of all, let us observe that the proof from~\cite{caro+yuster:00}
remains true for a much wider class of functions~$f$. Namely, a
function $f$ is called \emph{non-decreasing} if for any $m\le n$ we
have $f(m)\le f(n)$.

\bth{1} For any $n\ge 0$, $r\ge 3$ and non-decreasing $f:\I N\to\I R$
we have
 \beq{s2}
 \ex_f(n,K_r)=\ex_f'(n,K_r).
 \eeq
 \end{theorem}
 \bpf By~\req{s1} we have to prove the `$\le$'-inequality only. Let an
$F$-free graph $G$ achieve $\ex_f(n,F)$. By a theorem of Erd\H
os~\cite{erdos:70} there is an $(r-1)$-partite graph $H$ on the same
vertex set $V$ such that $d_H(x)\ge d_G(x)$ for every $x\in V$. We
have
 $$
 \ex_f(n,K_r)=\sum_{x\in V} f(d_G(x)) \le \sum_{x\in V} f(d_H(x)) \le
ex_f'(n,K_r),
 $$
 establishing the required.\epf

Caro and Yuster~\cite[Conjecture~6.1]{caro+yuster:00} posed the
problem of computing $\ex_{P_\mu}(n,F)$ for an arbitrary graph
$F$. Here we show that if $F$ is a fixed graph of chromatic number
$r\ge 3$ and $f$ is a `nice' function (including power functions
$P_\mu$ for $\mu\ge 0$) then in order to compute $\ex_f(n,F)$
asymptotically it is enough to consider complete $(r-1)$-partite
graphs only. Determining the optimal sizes of the $r-1$ parts may be a
difficult analytical task. (Bollob\'as and
Nikiforov~\cite{bollobas+nikiforov:04} investigate this problem for
the power function $P_\mu$.) But, combinatorially, the problem of
computing the asymptotics of $\ex_f(n,F)$ for such $f,F$ may be
considered as solved.

\section{Main Result}

Let us call a non-decreasing function $f:\I N\to\I R$
\emph{log-continuous} if for any $\e>0$ there is $\delta>0$ such that
for any $m,n\in \I N$ with $n\le m\le (1+\delta)n$ we have
 \beq{s3}
 f(m)\le (1+\e)f(n).
 \eeq
 For example, $P_\mu$ is log-continuous for any $\mu>0$ while the
exponent $x\mapsto \me^{x}$ is not.

\bth{2} Let $F$ be a fixed non-bipartite graph. Let $f:\I N\to\I R$
be an arbitrary non-decreasing and log-continous function. Then, as
$n\to\infty$,
 $$
 \ex_f(n,F)=(1+o(1))\, \ex_f'(n,F).
 $$
 \end{theorem}
 \bpf Let $r=\chi(F)\ge 3$, $\e>0$ be arbitrary, $n$ be large, and
$G$ achieve $\ex_f(n,F)$.

We will need the following result of Erd\H os, Frankl and
R\"odl~\cite[Theorem~1.5]{erdos+frankl+rodl:86}.

\begin{theorem}\label{th:efr86} For every $\e>0$ and a graph $F$,
there is a constant $n_0=n_0(\e,F)$ with the following property. Let
$G$ be a graph of order $n\ge n_0$ that does not contain $F$ as a
subgraph. Then $G$ contains a set $E'$ of less than $\e n^2$ edges
such that the subgraph $H=G-E'$ has no $K_r$, where
$r=\chi(F)$.\qed\end{theorem}

Theorem~\ref{th:efr86} is proved by applying Szemer\'edi's Regularity
Lemma so the lower bounds on $n_0$ are huge. Also, $\delta=\delta(\e)$
in~\req{s3} is implicit. Therefore, in what follows we make no attempt
to optimize the constants. In the proof we choose positive constants
$\gamma_1,\gamma_2,\gamma_3,\gamma_4$. The Reader can check that
indeed we can choose $\gamma_i$ sufficiently small (depending on
$F,f,\e,\gamma_1,\dots,\gamma_{i-1}$) so that all inequalities are
true for all large $n$.

First, let us observe that by the assumptions on $f$
 \beq{s4}
 \ex_f'(n,F)\ge e_f(T_{r-1}(n)) \ge n f(\floor{n(r-1)/r}) \ge \gamma_1 n f(n).
 \eeq

Let $V=V(G)$. Define $A=\{x\in V\mid d_G(x)\le \gamma_2n\}$ and
$B=V\setminus A$. The subgraph $G'=G[B]$ spanned by $B$ is of course
$F$-free. Theorem~\ref{th:efr86} gives us a $K_r$-free subgraph
$H'\subset G'$ with $e(G')-e(H')\le \gamma_4 n^2$. By the theorem of
Erd\H os~\cite{erdos:70} there is an $(r-1)$-partite graph on $B$
majorizing the degrees of $H'$. Extend it to a complete
$(r-1)$-partite graph on $V$ by arbitrarily splitting $A$ into $r-1$
almost equal parts.

The proof will be complete if we show that 
 \beq{aim}
 e_f(G)-e_f(H)\le \e e_f(G).
 \eeq

To prove~\req{aim} we estimate the contribution to $e_f(G)-e_f(H)$ by
various sets of vertices.

If, for example, $|A|\ge 2\gamma_2n$, then for any $x\in A$ we have
$d_H(x)\ge \floor{\frac{r-1}r\, |A|}\ge \gamma_2 n\ge d_G(x)$, ie.,\
the contribution of $A$ to the left-hand side of~\req{aim} is
non-positive. Otherwise, the contribution is at most
 $$
 \sum_{x\in A} f(d_G(x))\le |A| f(n)\le \frac{\e \gamma_1}4\, n
f(n)\le \frac{\e}4 e_f(G).
 $$
 (In the last inequality we used~\req{s4}.)

Let $C=\{x\in B\mid |\Gamma_G(x)\cap A|\ge \frac12|A|\}$. By counting
the degrees in $A$, we obtain $|C|\le 2\gamma_2n$ and thus the
$C$-contribution is also at most $\frac{\e}4\, e_f(G)$. Notice
that any vertex of $D=B\setminus C$ has at least as many $A$-neighbors
in $H$ as it has in $G$.

Define $K=\{x\in D\mid d_{H'}(x) \le d_{G'}(x)-\gamma_3 n\}$. Clearly, 
 $$
 |K|\le \frac{2(e(G')-e(H'))}{\gamma_3 n}\le
\frac{2\gamma_4}{\gamma_3}\, n.
 $$
 Again, $|K|$ is so small that $K$ contributes at most
 $\frac{\e}4\, e_f(G)$ to~\req{aim}.

Let us estimate the contribution of $K=D\setminus L$. For any vertex
$x$ of $D$ we have
 $$
 d_G(x)-d_H(x) \le d_{G'}(x)-d_{H'}(x) \le \gamma_3 n.
 $$ 
 All vertices in $K\subset B$ has $G$-degree at least $\gamma_2 x$. As
$\gamma_3$ is small compared with $\gamma_2$, we can assume by the
log-continuity of $f$ that $f(d_G(x))\le (1+\frac{\e}4)f(d_H(x))$ for
any $x\in K$. This completes the proof of~\req{aim} and the
theorem.\epf

\brm Taking $f:x\mapsto\log x$ we can also solve the problem of
maximizing $\prod_{x\in V(G)} d(x)$ over all $F$-free graphs $G$ of
order $n$. (However, please notice that the relative error here will
not be $1+o(1)$ but becomes such after taking the logarithm.) More
generally, we can maximize $\prod_{x\in V(G)} f(d(x))$ for any
non-decreasing $f$ such that $\log (f(x))$ is log-continuous; in
particular, this is true if $f$ itself is log-continous.\medskip

\section{Some Negative Examples}

In Theorem \ref{th:2} we do need some condition bounding the rate of
growth of $f$. For example, if $f$ grows so fast that $e_f(G)$ is
dominated by the contribution from the vertices of degree $n-1$, then
the conlusion of Theorem~\ref{th:2} is no longer true: for example,
for $K_3(2)$ (the blown-up $K_3$ where each vertex of $K_3$ is
duplicated) the value $\ex_f(n,K_3(2))=(3+o(1))\, f(n-1)$ cannot be
achieved by a bipartite graph.

In fact, one can get refuting examples of $f$ with moderate rate of
growth. For example, for any constant $c<1$ there is a non-decreasing
$f$ such that
 \beq{Df}
 \frac{f(n+1)}{f(n)} \le 1+n^{-c}
 \eeq
 for any $n$ and yet the conlusion of Theorem~\ref{th:2} does not hold
for this $f$. Let us demonstrate the above claim.

Let $c>0$. Choose $t$ such that for all large $n$ there is an
$K_{t,t}$-free graph $G_n$ of order $n$ with almost all vertices
having degree at least $n^c$ each. This $t$ exists by a construction
of Koll\'ar, R\'onyai and Szab\'o~\cite{kollar+ronyai+szabo:96}.

Let $F=K_3(2t-1)$ be a blown-up $K_3$. Take an arbitrary function $f$
satisfying~\req{Df} and the additional property that there is an
infinite sequence $n_1<n_2<\dots$ such that for any $k$ we have
 $$
 f(n_k+m_k)=f(n_k+m_k+1)=f(n_k+m_k+2)=\dots=f(2n_k),
 $$
 while $f(n_k)\le \frac12 f(n_k+m_k)$, where $m_k=\floor{\frac12
n^c}$. Such an $f$ exists: choose the numbers $n_k$ spaced far apart
(with $n_1$ being sufficiently large), let
$f(n+1)=f(n)$ except for $n_k\le n\le n_k+m_k$ we let $f(n+1)=2^{1/m_k}
f(n)$. Note that $2^{1/m_k}<1+\frac{1}{m_k}< 1 +n^{-c}$ so our $f$
does satisfy~\req{Df}.

On one hand, we have
 \beq{K33}
 \ex_f(2n_k,F)\ge (2+o(1))\, n_k f(n_k+m_k).
 \eeq
 Indeed, Let $G$ be obtained from the complete bipartite graph
$K_{n_k,n_k}$ by adding to each part the $K_{t,t}$-free graph
$G_{n_k}$ defined above. It is easy to see that $G\not\supset
F$. Almost all vertices of $G$ have degree at least $n_k+m_k$,
giving~\req{K33}. 

On the other hand, for any bipartite graph $H$ of order $2n_k$ at
least $n_k$ vertices will have degree at most $n_k$ and thus
 $$
 e_f(H)\le n_k f(2n_k-1) + n_k f(n_k)\le \frac32\, n_k
f(n_k+m_k).
 $$
 We obtain by~(\ref{eq:K33}) that $\ex_f(n,F)$ cannot always be
approximated by bipartite graphs.

\end{document}

\bibliography{oleg,general,graph,ex,sets}

\end{document}